\documentclass[11pt,twoside,center]{article}

\usepackage[dvips]{graphicx}

\title{\bf LIMIT POINTS AND HOPF BIFURCATION POINTS FOR A ONE - PARAMETER DYNAMICAL
SYSTEM ASSOCIATED TO THE LUO - RUDY I MODEL}

\author{C\u{a}t\u{a}lin Liviu Bichir\footnote{
{\tt catalinliviubichir@yahoo.com}, Rostirea Maths Research,
Regimentul 11 Siret 27, Gala\c{t}i, Romania} \and Adelina
Georgescu\footnote{Academy of Romanian Scientists, Splaiul
Independen\c{t}ei 54, Bucharest, Romania} \and Bogdan
Amuzescu\footnote{{\tt bogdan@biologie.kappa.ro}, Faculty of
Biology, University of Bucharest, Splaiul Independen\c{t}ei 91-95,
Bucharest, Romania. This research was partially supported from
grant PNCDI2 61-010 to M-LF by the Romanian Ministry of Education,
Research, and Innovation.} \and Gheorghe
Nistor\footnote{University of Pite\c{s}ti, Str. T\^{a}rgul din
Vale 1, Pite\c{s}ti, Romania} \and Marin
Popescu\footnote{University of Pite\c{s}ti, Str. T\^{a}rgul din
Vale 1, Pite\c{s}ti, Romania} \and Maria-Luiza
Flonta\footnote{Faculty of Biology, University of Bucharest,
Splaiul Independen\c{t}ei 91-95, Bucharest, Romania} \and
Alexandru Dan Corlan\footnote{Bucharest University Emergency
Hospital, Splaiul Independen\c{t}ei 169, Bucharest, Romania} \and
Istvan Svab\footnote{Faculty of Biology, University of Bucharest,
Splaiul Independen\c{t}ei 91-95, Bucharest, Romania}}
\date{~}
\begin{document}

\maketitle

\pagestyle{myheadings}
\markboth{Bichir, Georgescu, Amuzescu, Nistor et al.
}{LIMIT POINTS AND HOPF BIFURCATION POINTS}

\bigskip

\begin{abstract}
A one - parameter dynamical system is associated to the
mathematical problem governing the membrane excitability of a
ventricular cardiomyocyte, according to the Luo-Rudy I model. An
algorithm used to construct the equilibrium curve is presented.
Some test functions are used in order to locate limit points and
Hopf bifurcation points. Two extended systems allow to calculate
these points. The numerical results are presented in a bifurcation
diagram.
\newline
\end{abstract}

{\bf keywords:}  limit point, Hopf bifurcation point, Luo-Rudy I
model, arc-length-continuation method, Newton's method, computer
program.

\section{Introduction}
\label{BichirCL_sect_I} The present paper is one of a series of
research results, for the dynamical system associated to the
Luo-Rudy I model, obtained under the coordination of Acad. Adelina
Georgescu.

Mathematical models of cardiomyocyte electrophysiology based on
experimentally determined kinetics of ion currents encompass over
40 years, since the early attempts of Denis Noble
(\cite{BichirCL_Noble5}, \cite{BichirCL_Noble6}) to accurate
models of cell types in all regions of the heart that are now
being incorporated into anatomically detailed models of the whole
organ (\cite{BichirCL_Noble7}). These approaches introduced in
successive steps several time-dependent ion current components, as
well as a detailed dynamics of calcium in subcellular
compartments, comprising release and reuptake from the
sarcoplasmic reticulum (\cite{BichirCL_DiFrancescoNoble}),
specific calcium buffers (\cite{BichirCL_HilgemannNoble}), and the
electrogenic Na/Ca exchanger. The Luo-Rudy I model of ventricular
cardiomyocyte (\cite{BichirCL_LuoRudyI}), developed in the early
1990s starting from the Beeler-Reuter model
(\cite{BichirCL_BeelerReuter}), includes kinetics based on
single-channel recordings. All current components are described by
Hodgkin-Huxley type equations. This simplicity renders it adequate
for mathematical analysis using methods of linear stability and
bifurcation theory.
\newline
Nowadays, there exist numerous software packages for the numerical
study of finite - dimensional dynamical systems, for example
MATCONT, CL$_{-}$MATCONT, CL$_{-}$MATCONTM
(\cite{BichirCL_MATCONT}, \cite{BichirCL_Cl_MatContM}), AUTO
\cite{BichirCL_Doedel}. In this paper, numerical results are
obtained by using some new computer programs.

\section{Luo-Rudy I model}
\label{BichirCL_sect_II}

In spite of its simplicity, which comes from the fact that it does
not take into account earlier findings concerning
Ca\textsuperscript{2+} dynamics, the Luo-Rudy I model
\cite{BichirCL_LuoRudyI} proved to be very realistic,
incorporating data derived from single-channel recordings obtained
during the 1980s with the advent of the patch-clamp technique
 (\cite{BichirCL_Hamill}, \cite{BichirCL_Neher}). The model comprises
only three time and voltage-dependent ion currents (fast sodium
current, slow inward current, time-dependent potassium current)
plus three background currents (time-independent and plateau
potassium current, background current). The Luo-Rudy I model
reproduces a wealth of experimental findings, like: the fast
upstroke velocity of the action potential ($\dot{V}_{max}$ = 400
V/sec), the behavior of the rising phase, late repolarization
phase, and postrepolarization phase with changes in extracellular
potassium concentration $[K]_{o}$, monotonic Wenckebach patterns
and alternans at normal $[K]_{o}$, nonmonotonic Wenckebach
periodicities, aperiodic patterns, and enhanced supernormal
excitability resulting in unstable responses and chaotic activity
at low $[K]_{o}$.
\newline
Let us briefly describe how different experimental facts were
taken into account within the equations. The fast sodium current
($I_{Na}$) incorporates both a slow process of recovery from
inactivation and an adequate maximum conductance. The activation
(\textit{m}) and inactivation (\textit{h}) rates are adapted from
the Ebihara-Johnson $I_{Na}$ model based on data from chicken
embryo cardiac cells \cite{BichirCL_Ebihara}. Two inactivation
gates, fast and slow (\textit{h} and \textit{j}) were used to
render it compatible with single-channel data proving that near
threshold potentials sodium channels tend to open several times
during a depolarization (reopening phenomenon), and a significant
fraction of channels do not open by the time of peak inward
current. The start values of the slow inactivation gate \textit{j}
are obtained by setting $j_{\infty }=h_{\infty }$, as suggested by
Haas \textit{et al}. \cite{BichirCL_Haas}. The slow inward current
($I_{si}$) is represented exactly as in the Beeler-Reuter model.
The time-dependent potassium current ($I_{K}$) is controlled by a
time-dependent activation gate ($X$) and a time-independent
inactivation gate (${X}_{i}$) with inward rectification
properties, neither of which depends on $[K]_{o}$, while the
single-channel conductance is proportional to the square root of
$[K]_{o}$, as found in patch-clamp recordings on rabbit nodal
cells \cite{BichirCL_Shibasaki}. The time-independent potassium
current ($I_{K1}$) is different from $I_{K1}$ of the Beeler-Reuter
model, featuring two important properties discovered by Sakmann
and Trube using patch-clamp methods (\cite{BichirCL_Sakmann1},
\cite{BichirCL_Sakmann2}): a square-root dependence of
single-channel conductance on $[K]_{o}$, and a high selectivity
for potassium, as well as the inactivation gate K1 identified by
Kurachi in single-channel experiments \cite{BichirCL_Kurachi}.
Since this current inactivates completely during depolarization,
the model was supplemented with two other time-independent
potassium current components: a $[K]_{o}$ -insensitive plateau
current ($I_{Kp}$), simulating the single-channel properties of
the plateau current measured by Yue and Marban
\cite{BichirCL_Yue}, and a background current ($I_{b}$) with a
reversal potential $E_{b}$ = -59.87 mV. We should remark that,
although derived from single-channel experiments, the
conductances, gating and reversal potentials of these three
current components were adjusted using a parameter estimation
technique to fit the whole-cell time-independent potassium current
measured by Sakmann and Trube for different values of $[K]_{o}$.
\newline
The mathematical problem governing the membrane excitability of a
ventricular cardiomyocyte, according to the Luo-Rudy I model
(\cite{BichirCL_LuoRudyI}), is a Cauchy problem
\begin{equation}
\label{BichirCL_e01}
   u(0)=u_{0} \, ,
\end{equation}
for the system of first order ordinary differential equations
\begin{equation}
\label{BichirCL_e02}
   \frac{du}{dt}=\Phi(\eta,u) \, ,
\end{equation}
where $u$ $=$ $(u_{1},\ldots,u_{8})$ $=$ $(V$, $[Ca]_{i}$, $h$,
$j$, $m$, $d$, $f$, $X)$, $\eta$ $=$ $(\eta_{1},\ldots,\eta_{13})$
$=$ $(I_{st}$, $C_{m}$, $g_{Na}$, $g_{si}$, $g_{Kp}$, $g_{b}$,
$[Na]_{0}$, $[Na]_{i}$, $[K]_{0}$, $[K]_{i}$, $PR_{NaK}$, $E_{b}$,
$T)$, $M=R^{8}$, $\Phi:R^{13} \times M \rightarrow M$, $\Phi =
(\Phi_{1},\ldots,\Phi_{8})$,
\begin{eqnarray*}
   & \ & \Phi_{1}(\eta,u)=
      -\frac{1}{\eta_{2}}[I_{st}
         +\eta_{3}{u}_{3}{u}_{4}{u}_{5}^{3}({u}_{1}-E_{Na}(\eta_{7},\eta_{8},\eta_{13}))
         \\
   & \ & \qquad +\eta_{4}{u}_{6}{u}_{7}({u}_{1}-c_{1}+c_{2}\ln{u}_{2})
         \\
   & \ & \qquad +g_{K}(\eta_{10})X_{i}({u}_{1})
            ({u}_{1}-E_{K}(\eta_{7},\eta_{8},\eta_{9},\eta_{10},\eta_{11},\eta_{13}){u}_{8}
         \\
   & \ &  \qquad +g_{K1}(\eta_{10})K1_{\infty}(\eta_{9},\eta_{10},\eta_{13},{u}_{1})
             ({u}_{1}-E_{K1}(\eta_{9},\eta_{10},\eta_{13}))
         \label{BichirCL_e03} \\
   & \ & \qquad +\eta_{5}Kp({u}_{1})({u}_{1}-E_{Kp}(\eta_{9},\eta_{10},\eta_{13}))
          +\eta_{6}({u}_{1}-\eta_{12})] \, ,
         \\
   & \ & \Phi_{2}(\eta,u)=
      -c_{3}\eta_{4}{u}_{6}{u}_{7}({u}_{1}-c_{1}+c_{2}\ln{u}_{2})
      +c_{4}(c_{5}-{u}_{2}) \, ,
         \\
   & \ & \Phi_{\ell}(\eta,u)=
      \alpha_{\ell}({u}_{1})-(\alpha_{\ell}({u}_{1})+\beta_{\ell}({u}_{1}))u_{\ell} \, ,
   \ \ell=3,\ldots,8 \, .
\end{eqnarray*}
The definitions of variables $V$, $[Ca]_{i}$, $h$, $j$, $m$, $d$,
$f$, $X$, parameters $I_{st}$, $C_{m}$, $g_{Na}$, $g_{si}$,
$g_{K1}$, $g_{Kp}$, $g_{b}$, $[Na]_{0}$, $[Na]_{i}$, $[K]_{0}$,
$[K]_{i}$, $PR_{NaK}$, $E_{b}$, $T$, constants
$c_{1},\ldots,c_{5}$, functions $g_{K}$, $E_{Na}$, $E_{K}$,
$E_{K1}$, $E_{Kp}$, $K1_{\infty}$, $X_{i}$, $Kp$, $\alpha_{\ell}$,
$\beta_{\ell}$, default values of parameters and initial values of
variables in the Luo-Rudy I model are: $V$ - transmembrane
potential, $[Ca]_{i}$ - intracellular calcium concentration, $h$
and $j$ - fast and slow inactivation variable of $I_{Na}$
(probability of gate $h$ or $j$ to be open), $m$ - activation
variable of $I_{Na}$, $d$ and $f$ - activation and inactivation
variable of $I_{si}$, $X$ - activation variable of $I_{K}$,
$X_{i}$ - steady-state inactivation of $I_{K}$, $K1_{\infty}$ -
steady-state gating variable of $I_{K1}$, $Kp$ - steady-state
gating variable of $I_{Kp}$, $\alpha_{\ell}$ and $\beta_{\ell}$ -
voltage dependence of opening and closing rates expressed as
Boltzmann distribution functions for two distinct energy levels,
$I_{st}$ - steady depolarizing/hyperpolarizing applied current,
$C_{m}$ - membrane capacitance per unit area, $g_{Na}$ - maximal
conductance of fast voltage-gated sodium current (per unit area),
$g_{si}$ - maximal conductance of slow inward (calcium) current,
$g_{K}$ - maximal conductance of time-dependent potassium current
, $g_{K1}$ - maximal conductance of inward rectifier potassium
current, $g_{Kp}$ - maximal conductance of plateau potassium
current, $g_{b}$ - maximal conductance of background current,
$[Na]_{0}$, $[Na]_{i}$, $[K]_{0}$, $[K]_{i}$ - extra- and
intracellular concentrations of sodium and potassium, $PR_{NaK}$ -
sodium/potassium permeability ratio for $I_{K}$, $E_{Na}$,
$E_{K}$, $E_{K1}$, $E_{Kp}$, $E_{b}$ - reversal potentials of
$I_{Na}$, $I_{K}$, $I_{K1}$, $I_{Kp}$, $I_{b}$, $T$ - absolute
temperature.
\newline
For the continuity of the model, the reader is referred to
\cite{BichirCL_LivshitzRudy}, and for the treatment of the vector
field $\Phi$ singularities to \cite{BichirCL_CLB3}. $\Phi$ is of
class $C^{2}$ on the domain of interest.

\section{The one - parameter dynamical
system associated to the Luo - Rudy I model}
\label{BichirCL_sect_III}

We performed the study of the dynamical system associated with the
Cauchy problem (\ref{BichirCL_e01}), (\ref{BichirCL_e02}) by
considering only the parameter $\eta_{1}=I_{st}$ and fixing the
rest of parameters. Denote $\lambda=\eta_{1}=I_{st}$ and
$\eta_{\ast}$ the vector of the fixed values of
$\eta_{2},\ldots,\eta_{13}$. Let $F:R \times M \rightarrow M$,
$F(\lambda,u)=\Phi(\lambda,\eta_{\ast},u)$, $F =
(F_{1},\ldots,F_{8})$.
\newline
Consider the dynamical system associated with the Cauchy problem
(\ref{BichirCL_e01}), (\ref{BichirCL_e06}), where
\begin{equation}
\label{BichirCL_e06}
   \frac{du}{dt}=F(\lambda,u) \, .
\end{equation}
\newline
The equilibrium points of this problem are solutions of the
equation
\begin{equation}
\label{BichirCL_e07}
   F(\lambda,u)=0 \, .
\end{equation}
The existence of the solutions and the number were established by
graphical representation in \cite{BichirCL_CLB3}, for the domain
of interest. The equilibrium curve (the bifurcation diagram) was
obtained in \cite{BichirCL_CLB3}, via an arc-length-continuation
method (\cite{BichirCL_CLB10}) and Newton's method
(\cite{BichirCL_Gir_Rav1996}), starting from a solution obtained
by solving a nonlinear least-squares problem
(\cite{BichirCL_CLB10}) for a value of $\lambda$ for which the
system has one solution. In \cite{BichirCL_CLB3}, the results are
obtained by reducing (\ref{BichirCL_e07}) to a system of two
equations in $(u_{1},u_{2})$ $=$ $(V, [Ca]_{i})$. Here, we used
directly (\ref{BichirCL_e07}).

\section{Arc-length-continuation method and Newton's method for (\ref{BichirCL_e07})}
\label{BichirCL_sect_IV} Let us write the arc-length-continuation
method and Newton's method used to construct the equilibrium curve
of (\ref{BichirCL_e07}).
\newline
Glowinski (\cite{BichirCL_CLB10}, following H.B.Keller
\cite{BichirCL_Keller1}, \cite{BichirCL_Keller2}) chose a
continuation equation written in our case as
\begin{equation}
\label{BichirCL_e08}
   \sum\limits_{i=1}^{8}(\frac{du_{i} }{d s})^{2}
      + (\frac{d \lambda }{d s})^{2} = 1 \, ,
\end{equation}
where $s$ is the curvilinear abscissa.
\newline
Let $(\lambda_{\ast}, u_{\ast})$ be a solution of
(\ref{BichirCL_e07}) obtained by solving a nonlinear least-squares
problem (\cite{BichirCL_CLB10}) for a fixed value of $\lambda$
($\lambda=\lambda_{\ast}$) for which the system has one solution.
To solve (\ref{BichirCL_e07}), let us consider the extended system
formed by (\ref{BichirCL_e07}) and (\ref{BichirCL_e08}),
parameterized by $s$. Let $\triangle s$ be an arc-length step and
$u^{n} \cong u(n \triangle s)$. We have the algorithm (following
the case formulated in \cite{BichirCL_CLB10}): take
$\lambda(0)=\lambda^{0}=\lambda_{\ast}$, $u(0)=u^{0}=u_{\ast}$ and
suppose that $\frac{d \lambda (0)}{d s}$, $\frac{d u (0)}{d s}$
are given; for $n \geq 0$, assuming that $\lambda^{n-1}$,
$u^{n-1}$, $\lambda^{n}$, $u^{n}$ are known,
$(\lambda^{n+1},u^{n+1})$ is obtained by
\begin{equation}
\label{BichirCL_e09}
   F(\lambda^{n+1},u^{n+1})=0
\end{equation}
and
\begin{eqnarray}
   & \ & \sum\limits_{i=1}^{8}(u^{1}_{i}-u^{0}_{i})\frac{d u_{i} (0)}{d s}
      + (\lambda^{1}-\lambda^{0})\frac{d \lambda (0)}{d s}
      = \triangle s \ \textrm{if} \ n = 0 \, ,
         \nonumber \\
   & \ & \sum\limits_{i=1}^{8}(u^{n+1}_{i}-u^{n}_{i})
      \frac{u^{n}_{i}-u^{n-1}_{i}}{\triangle s}
      + (\lambda^{n+1}-\lambda^{n})\frac{\lambda^{n}-\lambda^{n-1}}{\triangle s}
         \label{BichirCL_e10} \\
   & \ & \qquad = \triangle s \, \ \textrm{if} \ n \geq 1 \, .
         \nonumber
\end{eqnarray}
\newline
In order to calculate $\frac{d \lambda (0)}{d s}$, $\frac{d u
(0)}{d s}$, we obtain the following relations.
From(\ref{BichirCL_e07}), we have
\begin{equation}
\label{BichirCL_e11}
   \sum\limits_{i=1}^{8}
      \frac{\partial F_{j}(\lambda^{0},u^{0})}{\partial u_{i}}
      \frac{d u_{i} (0)}{d s}
      + \frac{\partial F_{j}(\lambda^{0},u^{0})}{\partial \lambda}
      \frac{d \lambda (0)}{d s} = 0 \, , \, j=1,\ldots,8  \, .
\end{equation}
Let
\begin{equation}
\label{BichirCL_e12}
   \frac{d u_{i} (0)}{d s}=\hat{u}_{i}\frac{d \lambda (0)}{d s} \, , \, i=1,\ldots,8  \, .
\end{equation}
$\hat{u}$ $=$ $(\hat{u}_{1},\ldots,\hat{u}_{8})$ is the solution
of
\begin{equation}
\label{BichirCL_e13}
   \sum\limits_{i=1}^{8}
      \frac{\partial F_{j}(\lambda^{0},u^{0})}{\partial u_{i}}
      \; \hat{u}_{i} \;
      = \; - \; \frac{\partial F_{j}(\lambda^{0},u^{0})}{\partial \lambda}
      \, , \, j=1,\ldots,8  \, .
\end{equation}
From (\ref{BichirCL_e08}), we have
\begin{equation}
\label{BichirCL_e14}
   (\sum\limits_{i=1}^{8}\hat{u}_{i}+1)(\frac{d \lambda (0)}{d s})^{2}
   = 1 \, .
\end{equation}
\newline
In (\ref{BichirCL_e10}), let us denote $u^{\ast}=u^{0}$,
$\lambda^{\ast}=\lambda^{0}$, $u^{\ast\ast}=\frac{d u (0)}{d s}$,
$\lambda^{\ast\ast}=\frac{d \lambda (0)}{d s}$ if $n = 0$ and
$u^{\ast}=u^{n}$, $\lambda^{\ast}=\lambda^{n}$,
$u^{\ast\ast}=\frac{u^{n}-u^{n-1}}{\triangle s}$,
$\lambda^{\ast\ast}=\frac{\lambda^{n}-\lambda^{n-1}}{\triangle s}$
if $n \geq 1$.
\newline
The algorithm which we use to construct the branch of solutions
for (\ref{BichirCL_e07}) is the following: \\
1. given $\lambda^{0}$, $u^{0}$, solve (\ref{BichirCL_e13}) to
obtain
$\hat{u}$; \\
2. obtain $\frac{d \lambda (0)}{d s}$ from (\ref{BichirCL_e14})
and
$\frac{d u (0)}{d s}$ from  (\ref{BichirCL_e12}); \\
3. $\lambda(0)=\lambda^{0}=\lambda_{\ast}$ and
$u(0)=u^{0}=u_{\ast}$ are taken as discussed above; for $n \geq
0$, taking $\lambda^{n}$, $u^{n}$ as initial iteration, the
following algorithm based on Newton's method calculates
$\lambda^{n+1}$, $u^{n+1}$ : obtain $(\lambda^{m+1},u^{m+1})$
using
\begin{eqnarray}
   & \ & \sum\limits_{i=1}^{8}
      \frac{\partial F_{j}(\lambda^{m},u^{m})}{\partial u_{i}}
      \; u^{m+1}_{i}
      + \frac{\partial F_{j}(\lambda^{m},u^{m})}{\partial \lambda}
      \; \lambda^{m+1}
         \nonumber \\
   & \ & \qquad = \; \sum\limits_{i=1}^{8}
      \frac{\partial F_{j}(\lambda^{m},u^{m})}{\partial u_{i}}
      \; u^{m}_{i}
      + \frac{\partial F_{j}(\lambda^{m},u^{m})}{\partial \lambda}
      \; \lambda^{m}
         \nonumber \\
   & \ & \qquad  - F_{j}(\lambda^{m},u^{m}) \; , \, j=1,\ldots,8  \,
         \label{BichirCL_e15} \\
   & \ & \sum\limits_{i=1}^{8}u^{\ast\ast}_{i}u^{m+1}_{i}
      + \lambda^{\ast\ast}\lambda^{m+1}
      = \sum\limits_{i=1}^{8}u^{\ast}_{i}u^{\ast\ast}_{i}
      + \lambda^{\ast}\lambda^{\ast\ast} + \triangle s \, ;
         \nonumber
\end{eqnarray}
calculate the eigenvalues of the Jacobian matrix
$D_{u}F(\lambda^{n+1},u^{n+1})$ by the QR algorithm, calculate
$\psi_{LP}(\lambda^{n+1},u^{n+1})$ by (\ref{BichirCL_e16}), and
calculate
$\psi_{H}(\lambda^{n+1},u^{n+1})$ by (\ref{BichirCL_e17}); \\
4. the algorithm is stopped after an imposed number of iterations
for $n$.

\section{Limit points and Hopf bifurcation points}
\label{BichirCL_sect_V}

In order to locate the limit points and the Hopf bifurcation
points on the equilibrium curve of (\ref{BichirCL_e07}), two test
functions (\cite{BichirCL_Kuznetsov}, \cite{BichirCL_Govaerts},
\cite{BichirCL_Sey}, \cite{BichirCL_Sey2}), $\psi_{LP}$ and
$\psi_{H}$, are evaluated at each iteration
$(\lambda^{n+1},u^{n+1})$ of the algorithm from the end of section
\ref{BichirCL_sect_IV}, where
\begin{equation}
\label{BichirCL_e16}
   \psi_{LP}(\lambda,u)=det(D_{u}F(\lambda,u)) \, ,
\end{equation}
\begin{equation}
\label{BichirCL_e17}
   \psi_{H}(\lambda,u)=det(2 \, D_{u}F(\lambda,u) \,
      \odot \, I_{8}) \, .
\end{equation}
$\psi_{H}$ is evaluated using formula ($\ref{BichirCL_e18}$). For
a $n \times n$ matrix $A$ with elements $\{ a_{ij} \}$, the
following $m \times m$ matrix, $m=\frac{1}{2}n(n-1)$, is obtained
(\cite{BichirCL_Kuznetsov}, \cite{BichirCL_Govaerts}) based on the
definition of the bialternate product $\odot$ of two matrices,
\begin{equation}
\label{BichirCL_e18}
   (2 \, A \, \odot \, I_{n})_{(p,q),(r,s)} =
   \left\{ \begin{array}{l}
         -a_{ps} \, , \, \ \textrm{if} \ r=q \, ,
         \\
         a_{pr} \, , \, \ \textrm{if} \ r \neq p \ \textrm{and} \ s=q \, ,
         \\
         a_{pp}+a_{qq} \, , \, \ \textrm{if} \ r=p \ \textrm{and} \ s=q \, ,
         \\
         a_{qs} \, , \, \ \textrm{if} \ r=p \ \textrm{and} \ s \neq q \, ,
         \\
         -a_{qr} \, , \, \ \textrm{if} \ s=p \, ,
         \\
         0 \, , \, \ \textrm{otherwise} \, .
         \end{array} \right.
\end{equation}
The rows are labeled by the multi-index $(p,q)$ ($p=2,3,\ldots,n$,
$q=1,2,\ldots,p-1$), and the columns are labeled by the
multi-index $(r,s)$ ($r=2,3,\ldots,n$, $s=1,2,\ldots,r-1$).
\newline
If $\psi_{LP}$ has opposite signs at two points
$(\lambda^{n},u^{n})$ and $(\lambda^{n+1},u^{n+1})$, then a limit
point exists between $(\lambda^{n},u^{n})$ and
$(\lambda^{n+1},u^{n+1})$. If $\psi_{H}$ has opposite signs at
these two points, then it is possible that a Hopf bifurcation
points exists between them. The existence of a Hopf bifurcation
point is decided by studying the form of eigenvalues of
$D_{u}F(\lambda^{n},u^{n})$ and of
$D_{u}F(\lambda^{n+1},u^{n+1})$, since $\psi_{H}$ can be zero if
there is a pair of real eigenvalues of opposite sign and with
equal modulus. (We had in view the existence of a pair of complex
conjugate eigenvalues, for each matrix, and a change of the sign
of their real part).
\newline
For the cases where the test functions detect a limit point or a
Hopf bifurcation point, we retain the results of one of the
iterations $(\lambda^{n},u^{n})$ or $(\lambda^{n+1},u^{n+1})$ of
the algorithm presented at the end of section
\ref{BichirCL_sect_IV}, namely the iteration where the modulus of
the test function is smaller. These results are an initial
iteration for Newton's method applied to one of the equations
(\cite{BichirCL_Sey1}, \cite{BichirCL_Sey}, \cite{BichirCL_Sey2})
\begin{equation}
\label{BichirCL_e19}
   G(\lambda,u,h)=0 \, ,
\end{equation}
\begin{equation}
\label{BichirCL_e20}
   H(\lambda,\beta,u,h,g)=0 \, .
\end{equation}
\newline
The components $(\lambda,u)$ of the solution of
(\ref{BichirCL_e19}) represent a limit point of
(\ref{BichirCL_e07}). The components $(\lambda,u)$ of the solution
of (\ref{BichirCL_e20}) represent a Hopf bifurcation point of
(\ref{BichirCL_e06}). $H$ and $G$ are defined by
\begin{equation}
\label{BichirCL_sist_LP}
   G:R^{2 \cdot 8+1} \rightarrow R^{2 \cdot 8+1} , \
   G(\lambda,u,h)=\left[\begin{array}{l}
      F(\lambda,u) \\
      D_{u}F(\lambda,u)h\\
      h_{k} - 1
      \end{array}\right] \, ,
\end{equation}

\begin{equation}
\label{BichirCL_sist_H}
   H:R^{3\cdot 8+2} \rightarrow R^{3\cdot 8+2} , \
   H(\lambda,\beta,u,h,g)=\left[\begin{array}{l}
      F(\lambda,u) \\
      D_{u}F(\lambda,u)h+\beta g\\
      D_{u}F(\lambda,u)g-\beta h\\
      h_{k} - 1\\
      g_{k}
      \end{array}\right] \, ,
\end{equation}
where $k$ is a fixed index, $1 \leq k \leq 8$. The extended system
(\ref{BichirCL_e20}) determines a Hopf bifurcation point
$(\lambda,u)$, a pair of purely imaginary eigenvalues $\pm \beta
i$ of the Jacobian matrix in $(\lambda,u)$, and a nonzero complex
vector $h+ig$.
\newline
Newton's method applied to the equation (\ref{BichirCL_e19}) is:
let $v^{0}$ $=$ $(\lambda^{0}$, $u^{0}$, $h^{0})$  be an initial
iteration, where $\lambda^{0}$, $u^{0}$ are as specified above and
$h^{0}=(1,0,0,0,0,0,0,0)$; $k=1$; for $m \geq 0$, calculate
$v^{m+1}$ $=$ $(\lambda^{m+1}$, $u^{m+1}$, $h^{m+1})$ using
\begin{equation}
\label{BichirCL_e21}
   DG(v^{m})(v^{m+1}-v^{m})
       \, = \, - \, G(v^{m}) \, .
\end{equation}
\newline
Newton's method applied to the equation (\ref{BichirCL_e20}) leads
to: let $w^{0}$ $=$ $(\lambda^{0}$, $\beta^{0}$, $u^{0}$, $h^{0}$,
$g^{0})$ be an initial iteration, where $\lambda^{0}$, $u^{0}$ are
as specified above, $\beta^{0}$ is the positive complex part of
one of the two complex conjugate eigenvalues of the Jacobian
matrix in $(\lambda^{0}$, $u^{0})$, $h^{0}=(1,0,0,0,0,0,0,0)$,
$g^{0}=(0,0,0,0,0,0,0,0)$; $k=1$; for $m \geq 0$, calculate
$w^{m+1}$ $=$ $(\lambda^{m+1}$, $\beta^{m+1}$, $u^{m+1}$,
$h^{m+1}$, $g^{m+1})$ using
\begin{equation}
\label{BichirCL_e22}
   DH(w^{m})(w^{m+1}-w^{m})
       \, = \, - \, H(w^{m}) \, .
\end{equation}

\begin{figure}
   \raggedleft
   \includegraphics[height=4cm,width=6cm]{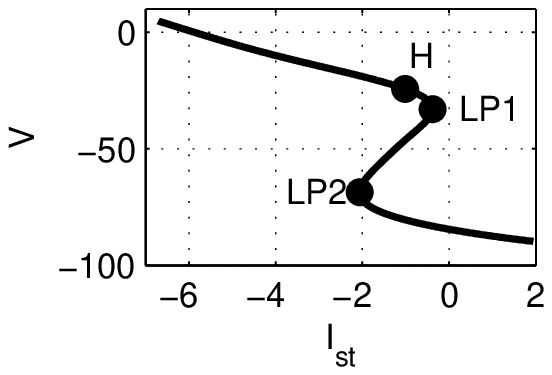}
   \centering
   \includegraphics[height=4cm,width=6cm]{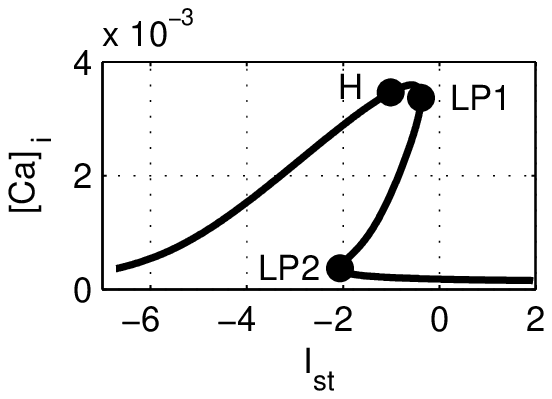}
   \raggedleft
   \includegraphics[height=4cm,width=6cm]{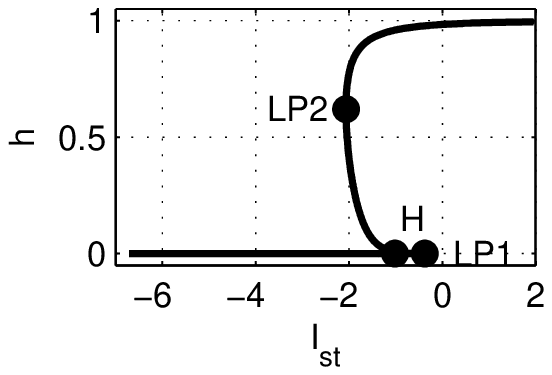}
   \centering
   \includegraphics[height=4cm,width=6cm]{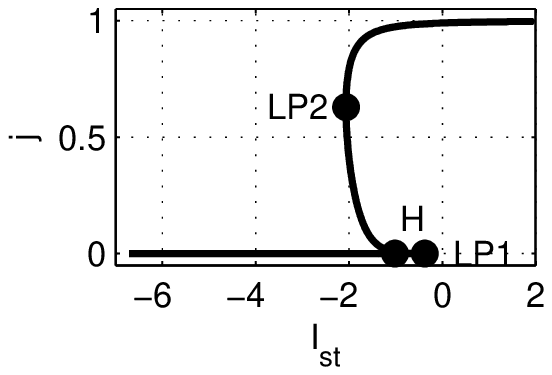}
   \raggedleft
   \includegraphics[height=4cm,width=6cm]{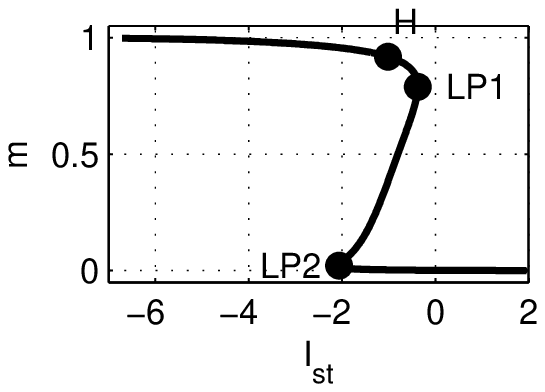}
   \centering
   \includegraphics[height=4cm,width=6cm]{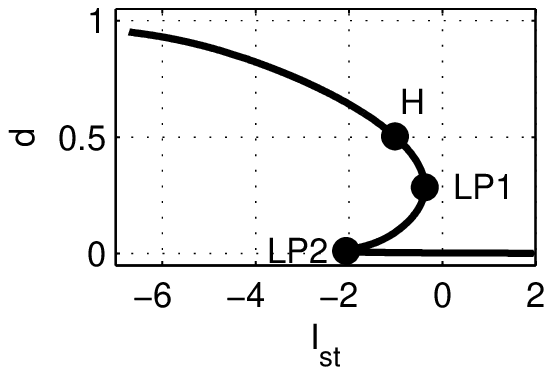}
   \raggedleft
   \includegraphics[height=4cm,width=6cm]{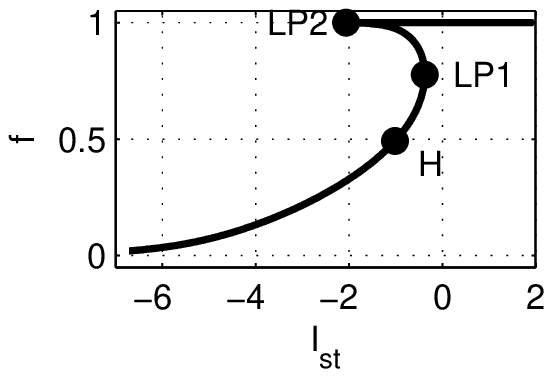}
   \centering
   \includegraphics[height=4cm,width=6cm]{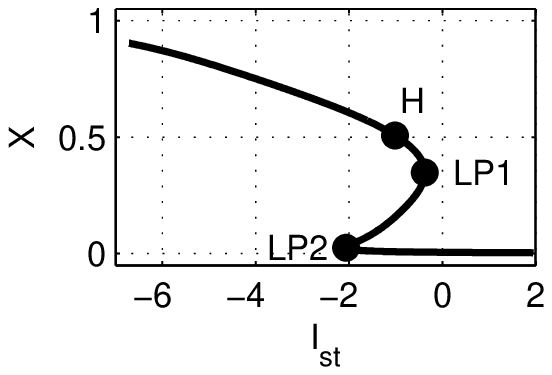}
   \raggedleft
   \includegraphics[height=4cm,width=6cm]{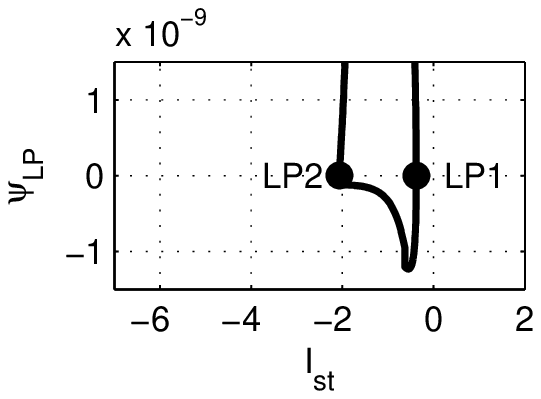}
   \centering
   \includegraphics[height=4cm,width=6cm]{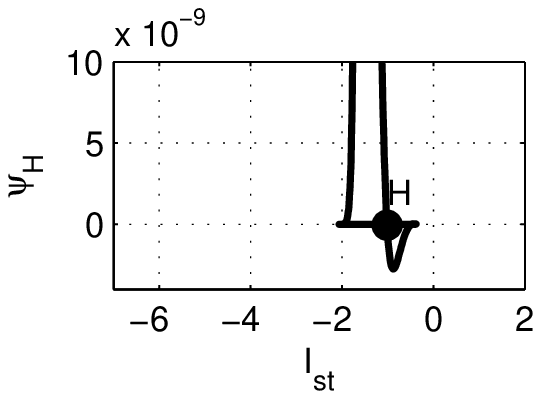}
   \caption{The bifurcation diagram and the variation of $\psi_{LP}$ and
$\psi_{H}$. There exist two limit points (LP1, LP2) and a Hopf
bifurcation point (H).}
   \label{BichirCL_diagr_bif}
\end{figure}

\section{Numerical results}
\label{BichirCL_sect_VI}

Based on the computer program for \cite{BichirCL_CLB2}, the
algorithm at the end of section \ref{BichirCL_sect_IV} was
transformed into a new computer program. Two new computer programs
were written in order to solve (\ref{BichirCL_e19}) and
(\ref{BichirCL_e20}) by (\ref{BichirCL_e21}) and
(\ref{BichirCL_e22}) respectively.
\newline
We took $C_{m}=1$, $g_{Na}=23$, $g_{si}=0.09$, $g_{K}=0.282$,
$g_{K1}=0.6047$, $g_{Kp}=0.0183$, $G_{b}=0.03921$, $[Na]_{0}=140$,
$[Na]_{i}=18$, $[K]_{0}=5.4$, $[K]_{i}=145$, $PR_{NaK}=0.01833$,
$E_{b}=-59.87$, $T=310$.
\newline
The equilibrium curves (the bifurcation diagram) are presented in
figure \ref{BichirCL_diagr_bif} for the domain of interest. The
variation of the functions $\psi_{LP}$ and $\psi_{H}$ are also
represented. "LP1" and "LP2" indicate two limit points. "H"
indicates a Hopf bifurcation point.
\newline
The behavior of the dynamical system changes with the value of the
variable parameter $I_{st}$. The two limit points separate three
branches of stationary solutions (first graph in Fig. 1). While
solutions on the middle branch (according to values of $V$) are
always unstable and those on the lower branch are always stable,
on the upper branch there is a region where the system features
oscillatory behavior. Oscillations are either damped, at the left
of the Hopf bifurcation point, or amplified until the system falls
on the lower branch of solutions to the right of the Hopf
bifurcation. Amplified oscillations represent early
afterdepolarizations (EADs), a condition prone to result in
life-threatening arrhythmias. A recent study of the Luo-Rudy I
dynamical system for variable relaxation time constants of the
gating variables $d$, $f$, and $X$ (\cite{BichirCL_Tran}), has
proved that oscillations resulting in EADs appear above a Hopf
bifurcation point for a fast subsystem, comprising the variables
$V$, $d$, and $f$. Moreover, these EADs can result in chaotic
behavior when the system is paced at a constant cycle length. The
same group has shown on detailed three-dimensional ventricular
electrophysiology models that EADs occurring in certain regions
can synchronize, resulting in polymorphic ventricular tachycardia
or torsades-de-pointes (\cite{BichirCL_Sato}).
\newline
In conclusion, our study, focused on analysis of the Luo-Rudy I
system as a whole in conditions of variable parameter $I_{st}$,
has identified two limit points and a Hopf bifurcation point,
separating different regions of stability, some of them featuring
amplified self-sustained oscillations defined as EADs on the time
trajectories, and which may result in dangerous ventricular
arrhythmias by synchronization. In contrast to the majority of
previous arrhythmogenesis studies, which attributed the generation
of this phenomenon to an individual condition, such as altered
gating of an ion channel type due to gene mutations or modulation
by physiological or pharmacological mechanisms, our results prove
that arrhythmogenesis may result as an emergent feature of the
system as a whole, and not of its individual components.

\enddocument